\documentclass[12pt]{article}

\usepackage{style}

\date{}

\title{Generalized and Convex Poincaré Inequalities via Transport-Divergence Inequalities}

\author[*]{Thibaut Le Gouic}
\author[*]{Magali Tournus}
\affil[*]{Aix-Marseille Univ., CNRS, I2M, UMR7373, Centrale Marseille, 13451 Marseille, France}

\begin{document}
\maketitle

\begin{abstract}
We investigate the relationship between generalized $(q',p')$-Poincar\'e inequalities and transport-divergence inequalities involving $\chi^q$-type divergences. Extending the classical correspondence between the $(2,2)$-Poincar\'e inequality and quadratic transport-$\chi^2$ inequalities, we prove that a $(q',p')$-Poincar\'e inequality is equivalent to a transport-divergence inequality 
for conjugates $(p,p')$ and $(q,q')$.
We then study the convex counterpart of these inequalities. Using weak barycentric optimal transport and convex order, we establish an equivalence between convex $(q',p')$-Poincar\'e inequalities and weak transport-divergence inequalities involving a projected $\chi_q$ divergence. We further characterize the minimizers arising in this projection problem through a geometric condition on the relative density.
In dimension one, we analyze the geometric consequences of $(q',p')$-Poincar\'e inequalities by identifying sharp conditions on the behavior of the density both at infinity and near isolated zeros.
 We show that when $q'>p'$, the concentration behavior becomes polynomial rather than exponential at infinity.
Finally, we interpret these inequalities through the lens of optimization in Wasserstein space. In particular, the equivalence between transport-divergence and Poincar\'e inequalities appears as an analogue of the equivalence between quadratic growth and Polyak--\L ojasiewicz inequalities in finite-dimensional optimization.
\end{abstract}

\begin{refsection}
\section{Introduction}
Since the seminal works \cite{OV2000}, deep connections have been established between classical functional inequalities, notably the Talagrand inequality and the log-Sobolev inequality. These inequalities are tightly linked: the log-Sobolev inequality implies Talagrand’s inequality, while the converse holds under log-concavity assumptions. Such relationships play a central role in the study of measure concentration, in particular in deriving Gaussian concentration bounds and convergence rates for associated Markov semigroups.

This perspective has been further clarified through the development of the theory of metric measure spaces with Ricci curvature bounded from below, due to Sturm, Lott, and Villani \cite{Sturm06,LV09}, as well as through Otto’s calculus on the Wasserstein space. Within this framework, these functional inequalities admit natural analogues in $\R^d$, revealing an underlying geometric structure (see Section~\ref{sec:optim}).

More recently, \cite{ding2015note,ledoux2018remarks} and subsequently \cite{liu2020}, established similar connections between transport–$\chi^2$ inequalities and the Poincaré inequality. These inequalities correspond to weaker concentration properties—typically exponential rather than Gaussian—and therefore hold in more general settings.

In this work, we propose to extend these results in two directions. First, we investigate inequalities associated with even weaker concentration regimes, namely polynomial. Second, we restrict attention to more structured classes of functions for the Poincaré inequality, namely convex functions. This leads us to study convex Poincaré inequalities as well as the $(q',p')$-Poincaré inequalities introduced in particular by Emanuel Milman \cite{milmanRoleConvexityIsoperimetry2009}, and to explore their connections with suitable transport inequalities.

\subsection{Poincaré inequality and transportation variance inequality}

Poincaré inequalities are fundamental tools in analysis, expressing how the behavior of a function can be controlled by its variations. 
In the classical setting, it writes
\begin{equation*}
    C_{2,2}^P \|f\|_{L^2(\mu)} \leq  \|\nabla f\|_{L^2(\mu)}, \qquad \forall f\text{ s.t. }\dst\int f\ud\mu=0,
\end{equation*}

They appear in a wide range of mathematical contexts and have numerous applications.
The Poincaré constant $C^P_{2,2}$ coincides with the (square root of the) spectral gap (that is, the first non-trivial eigenvalue) of the Fokker--Planck operator associated with the invariant measure $\mu$.
The spectral gap governs the rate of exponential convergence in $L^2(\mu)$ toward equilibrium of solutions to the Fokker--Planck equation $
\partial_t \mu_t=\nabla \cdot (\mu_t \nabla V)+\sqrt{2}\,\Delta \mu_t,$
 where $\mu = e^{-V}$.
A larger value of $C^P_{2,2}$ corresponds to faster convergence to equilibrium.
From a sampling perspective, when approximating $\mu$ using the Langevin diffusion—the particle-level counterpart of the Fokker--Planck equation—a small Poincaré constant typically signals the presence of geometric or energetic bottlenecks. These bottlenecks hinder efficient exploration of the support of $\mu$ by the Langevin dynamics and lead to slow mixing.

The Poincaré inequality is tightly linked to measure concentration phenomena \cite{Ledoux2001_book}. 
On $\R^d$, a Poincaré inequality for a measure $\mu$ implies exponential concentration of the measure 
and conversely exponential concentration yields a Poincaré inequality \cite{GRS11}. This correspondence extends to tensorized measures.
A measure $\mu$ is said to satisfy a dimension-free concentration property if its tensor product $\mu^{\otimes n}$ satisfies a concentration inequality with constants independent of $n$. It has been known since 1983 that any probability measure $\mu$ satisfying a Poincaré inequality also satisfies a dimension-free exponential concentration property \cite{gromov1983topological}.
In 2015, the reverse implication was established: if $\mu$ satisfies any form of dimension-free concentration property (Gaussian, exponential, or polynomial), then $\mu$ satisfies a Poincaré inequality \cite{GRS15}. The authors of \cite{GRS15} further extended this result by showing that dimension-free exponential concentration restricted to convex sets is equivalent to a Poincaré inequality restricted to convex functions.

 In 2012, \cite{jourdain2012equivalence} proved the
equivalence in dimension $n=1$ between  the classical $(2,2)$-Poincaré inequality and the transportation-variance
 inequality
 $$ C_{2,2}^T W_2(\mu,\nu) \leq \sqrt{\chi^2(\nu|\mu)}, \quad \nu\in\mathcal{P}_2(\R).$$
In 2015, \cite{ding2015note} proved the implication from the Poincaré inequality to the transport-variance inequality in arbitrary dimension, while Ledoux \cite{ledoux2018remarks} proved the reverse implication.
However, the proof of \cite{ding2015note} is somewhat technical and the constants are not optimal.
In 2020, the result was improved in \cite{liu2020} yielding the sharper constant $C_{2,2}^P = \sqrt{2} C_{2,2}^T$.
Those earlier proofs are based on the gradient flow structure in Wasserstein space $W_2$, which does not extend naturally beyond the exponent $p=2$. In contrast, our proofs rely on a duality argument and are more direct, yielding sharp constants and extending to general exponents.

The motivation to prove the equivalence between Poincaré inequality and Transportation-Variance inequality partially arises from the problem of characterizing the exponential decay of the quadratic Wasserstein distance along the heat flow: the exponential decay of the Fokker--Planck semigroup in $L^2(\mu)$-norm (i.e.\ in $\chi^2$ divergence) is equivalent to the $(2,2)$-Poincaré inequality. Instead of attacking directly the analogous question in Wasserstein space---namely, finding an equivalent inequality characterizing exponential decay of $W_2$---one natural approach is to compare $W_2$ with the square root of the variance.
The transport-variance inequality is then directly used in \cite{lehec2023langevin}
to prove the convergence in $W_2$ of the Langevin algorithm
assuming a Poincaré inequality for the target measure $\mu$.

\subsection{\texorpdfstring{($q',p'$)}{(q',p')}-Poincaré inequality}

Classical Poincaré inequalities impose an exponential decay of the measure and are therefore restrictive. In many situations of interest, the measure $\mu$ has heavier tails, and weaker functional inequalities are needed to capture its concentration properties. This motivates the study of $(q',p')$-Poincaré inequalities, which correspond to a broader range of tail behaviors.

\begin{definition}\label{def:poincaré}
A measure $\mu$ satisfies a \emph{$(q',p')$-Poincaré} inequality if there exists a constant $C>0$ such that for all $f \in L^{q'}(\mu)\cap C^1(\R^d)$ such that $\int |f|^{q'-2}f\ud\mu=0$
\begin{equation}\label{eq:Poincaré}
C\|f\|_{L^{q'}(\mu)}\le \|\nabla f\|_{L^{p'}(\mu)}.
\end{equation}
The largest such $C$ is denoted by $C^P_{q',p'}$.
\end{definition}

The condition $\int |f|^{q'-2}f\ud\mu=0$
is natural, as it is equivalent to translating $f$ by a constant $c$ chosen so as to minimize
$\|f-c\|_{L^{q'}(\mu)}$.
Observe that when $q'=2$, it reduces to $\esp_\mu[f]=0$, while for $q'=1$, it corresponds to requiring that $0$ is a median of $f$. 
Due to monotonicity of $L^r(\mu)$-norms with $r$, as soon as a measure $\mu$ satisfies a $(q',p')$-Poincaré inequality, then it satisfies a $(\tilde q, \tilde p)$-Poincaré inequality for any $\tilde q \leq q'$ and any $\tilde p \geq p'$.
In other words, for any measure, $C^P_{q',p'}\le C^P_{\tilde q,\tilde p}$ whenever $\tilde q \leq q'$ and $\tilde p \geq p'$.
In a breakthrough, \cite{milmanRoleConvexityIsoperimetry2009} proved that for log-concave measures $\mu$, all $(q',p')$-Poincaré inequalities are equivalent in the range $1 \leq q' \leq p'$: for any other $1 \leq q' \leq p'$, we have  $C_{q,p}^P \leq C q' C_{q',p'}^P$ for some universal constant $C$.
In \cite{cattiaux2020poincare}, the stability of Poincar\'e inequalities under perturbations of the underlying measure is investigated. In particular, several distances between probability measures are considered, including the total variation, Wasserstein, and bounded-Lipschitz distances. The authors show that if the distance between $\mu$ and $\nu$ is sufficiently small, then the Poincar\'e constant of $\mu$ controls that of $\nu$.
It is also proved that if the convolution $\nu \star \gamma$ satisfies a Poincar\'e inequality for some Gaussian measure $\gamma$, then $\nu$ itself satisfies a Poincar\'e inequality. This observation justifies restricting the study of Poincar\'e inequalities to smooth measures.

We note that if the support of the measure $\mu$ is not connected, then all Poincaré constants $C^P_{q',p'}$ vanish.
In particular, finite sums of Dirac masses fail to satisfy any Poincaré inequality.
However, a measure whose density vanishes at an isolated point may still satisfy a Poincaré inequality, provided that the density does not degenerate too rapidly near that point. As we shall see in Section~\ref{sec:role}, the admissible rate of vanishing depends on the parameter $p'$: larger values of $p'$ allow for stronger degeneracies of the density near the singular point.

These local behaviors,  which classical Poincar\'e inequalities are sensitive to, become invisible when one restricts attention to convex test functions.
This partly justifies why convex Poincaré inequalities have been studied (see \cite{Talagrand90,Sam03}).
\begin{definition}[Convex $(q',p')$-Poincaré inequality]
A measure $\mu$ is said to satisfy a \emph{convex $(q',p')$-Poincar\'e inequality} if there exists a constant $C>0$ such that for every convex function $f \in L^{q'}(\mu)\cap C^1(\R^d)$ satisfying
\[
\int |f|^{q'-2}f\,\ud\mu=0,
\]
one has
\[
C\|f\|_{L^{q'}(\mu)}\le \|\nabla f\|_{L^{p'}(\mu)}.
\]
The largest such constant is denoted by $C^{cP}_{q',p'}$.
\end{definition}

Convex Poincar\'e inequalities are weaker than their classical counterparts; in particular, every measure with bounded support satisfies a convex Poincar\'e inequality.
Notably, the behavior at infinity of a measure satisfying a convex $(q',p')$-Poincar\'e inequality is the same as in the classical $(q',p')$-Poincar\'e setting.


\subsection{Transport-divergence inequalities}

We recall the definition of the $p$-Wasserstein distance between $\mu$ and $\nu$ in $\mathcal{P}_p(\R^d)$:
\begin{equation*}
W_p^p(\mu,\nu)=\inf_{\pi\in\Pi(\mu,\nu)} \dst\int_{\R^d\times \R^d} |x-y|^p \ud\pi(x,y),
\end{equation*}
where $\Pi(\mu,\nu)$ is the set of couplings between $\mu$ and $\nu$.

We also define the $\chi^q$ divergence, in analogy with the classical $\chi^2$ divergence in statistics, by replacing the square by a general exponent $q\ge 1$:
\begin{equation*}
\chi^q(\nu|\mu )=
\begin{cases}
    &\left\|\dfrac{\nu}{\mu}-1\right\|_{L^q(\mu)}^q \quad  \text{ if } \nu\ll \mu,\\
    & +\infty \quad \qquad \qquad  \text{ otherwise.}
\end{cases}
\end{equation*}
\begin{definition}[Transport-divergence inequality]
A measure $\mu$ satisfies a \emph{$(p,q)$-transport divergence} inequality if there exists $C>0$ such that for all $\nu\in\mathcal{P}_p(\R^d)$ absolutely continuous w.r.t. $\mu$,
\[
CW_p(\mu,\nu)\le\left[\chi^q(\nu|\mu )\right]^{1/q}.
\]
The largest such $C$ is denoted by $C^T_{p,q}$.
\end{definition}

The inequality for the case $p=1,q=2$ has been introduced in \cite{ding2014wasserstein} to study polynomial concentration inequalities, with $\alpha$-Rényi divergence on the right hand side, which coincides with $\chi^2$ when $\alpha=2$.
They proved in their Theorem 1.6 that such transport-divergence inequality holds with $p=1$ and $q=2$ when $\mu$ admits a moment of order $4$.
In \cite[Theorem 1.1]{bobkov2015optimal}, they improved the result showing that such inequality holds if and only if $\mu$ has a finite moment of order $2$, i.e. $\mu\in\mathcal{P}_2(\R^d)$.

Additionally, we investigate convex Poincaré inequalities, where the test functions $f$  are convex.
The corresponding transport-divergence inequalities involve weak optimal transport. 
\begin{definition}[Weak barycentric Optimal Transport]\label{def:wbo}
The weak barycentric $p$-optimal transport between
$\mu$ and $\nu$ is defined as
\begin{equation*}
\overline{\mathcal{T}}_p(\nu|\mu)= \inf_{\pi \in \Pi(\mu,\nu)} \dst\int_{\R^d} \left| \dst\int_{\R^d} y \ud\eta_x(y) -x \right|^p \ud\mu(x),
\end{equation*}
where $\eta$ is the disintegration kernel such that
\begin{equation*}
\ud\pi(x,y) = \ud\mu(x)\, \ud\eta_x(y).
\end{equation*}
\end{definition}
The bar in the notation stands for {\it bar}ycentric.

Convex Poincaré inequalities and weak transport-divergence inequality have been explored, for instance in \cite{gozlan2018characterization}
 where the authors prove the equivalence of convex Poincaré inequalities with a weak transport-entropy inequality on the real line, in the case $p=2$ and for an appropriate local cost. 
The result was then extended in arbitrary dimension $n$ in \cite{adamczak2019convex}.

\subsection{Main results}
As stated above, our goal is to extend the result of \cite{ding2015note,ledoux2018remarks,jourdain2012equivalence,liu2020} on the equivalence between $(2,2)$-Poincaré inequalities and transport-$\chi^2$-divergence inequalities to $(q',p')$-Poincaré inequalities and transport-divergence inequalities and their convex counter parts.
In order to introduce our main results, we introduce $p,p',q,q'> 1$ satisfying 
\begin{align*}
\frac{1}{p}+\frac{1}{p'}=1 \quad \text{ and }\quad  \frac{1}{q}+\frac{1}{q'}=1.
\end{align*}

Recall that for a measure $\mu$, the best possible constant for the $(q',p')$-Poincaré inequality to hold is $C^P_{q',p'}$ and $C^T_{p,q}$ is the best possible constant for the transport divergence inequality to hold.
Thus we prove the following.
\begin{theorem}[Poincaré and Transport Divergence inequalities]
\label{thm:Poincare_transport_intro}
\begin{equation*}
 C^T_{p,q} \leq C^P_{q',p'} \leq p C^T_{p,q}.
\end{equation*}
Furthermore, the bounds have an optimal dependence on $p$, as we exhibit examples showing that no sharper estimates are possible.
\end{theorem}

Thanks to the tools we develop in the proofs, we are able to prove the same results for the convex Poincaré inequalities and weak transport divergence inequalities, in Theorem~\ref{thm:convex_Poincare_to_transport} and Theorem~\ref{thm:convex_transport_to_Poincare}.

The paper is organized as follows. Section \ref{sec:proof} is devoted to the proof of Theorem \ref{thm:Poincare_transport_intro}, while Section \ref{sec:convex} contains the proof of its convex equivalent.
In Section \ref{sec:role}, we emphasize the role played by the parameters $q'$ and $p'$ in
connection with concentration properties of the measure and local behavior near zeros in dimension one.
Finally, in Section \ref{sec:optim}, functional inequalities such as Poincaré and logarithmic Sobolev inequalities are interpreted through the lens of classical optimization inequalities in $\R^d$
. From this perspective, the hierarchy between functional inequalities can be viewed as an analogue of the standard hierarchy between classical optimization inequalities.

\section{Poincaré inequality and transport-divergence inequality}
\label{sec:proof}

This section is devoted to proving the two implications between $(q',p')$-Poincaré inequalities and transport-divergence inequalities.
The proofs are based on dual forms of each terms of the inequalities.
To show the first implication, we apply the Poincaré inequality to the Kantorovitch potential that appear in the dual form of the Wasserstein distance, and recover a transport-divergence inequality.
For the second implication, we apply a transport-divergence inequality to a measure whose relative density is given by the function that appears in the dual form of $L^q$-norms.
The fact that we need to apply the inequalities to different quantities when proving one implication and the other leads to a factor $p$ between the two results, which turns out to be of optimal order as we will see in the Remarks~\ref{rmk:counterT} and \ref{rmk:counterP} below.

\subsection{Poincaré implies transport divergence}

Here we show that for $p,p',q,q'>1$ 
\[
C^P_{q',p'}\le pC^T_{p,q}.
\]
\begin{theorem}
\label{thm:Poincare_transport}
If for all $f\in L^{q'}(\mu)\cap C^1(\R^d)$ such that $\int |f|^{q'-2}f\ud\mu=0$,
\[
C\|f\|_{L^{q'}(\mu)}\le\|\nabla f\|_{L^{p'}(\mu)},
\] then for all $\nu\in\mathcal{P}_p(\R^d)$,
    \[
    \frac{C}{p}W_p(\mu,\nu)\le \left[\chi^q(\nu|\mu )\right]^{1/q}.
    \]
\end{theorem}

\begin{proof}
We consider $\nu \in \mathcal{P}(\R^d)$ such that $\nu$ is absolutely continuous with respect to $\mu$. We use the duality formulation to write
\begin{equation*}
    W_p^p(\mu,\nu) = \dst\int \varphi^c \ud\nu + \dst\int \varphi \ud\mu,
\end{equation*}
where $\varphi(x) + \varphi^c(y)\leq
\|x-y\|^p$.
We notice that this inequality implies that for all $x\in \supp(\mu)$, $\varphi^c(x) \leq - \varphi(x)$, and then, since $\supp(\nu) \subset \supp(\mu)$,
\begin{equation*}
\begin{aligned}
    W_p^p(\mu,\nu) \leq \dst\int \varphi \ud(\mu-\nu)  &=\dst\int \varphi\left(1- \frac{\nu}{\mu} \right) \ud\mu\\
    & \leq \|\varphi\|_{L^{q'}(\mu)}\left\|\dfrac{\nu}{\mu}-1\right\|_{L^q(\mu)}  \\
    &\leq \dfrac{1}{C^P_{q',p'}} \|\nabla\varphi\|_{L^{p'}(\mu)}\left\|\dfrac{\nu}{\mu}-1\right\|_{L^q(\mu)},
    \end{aligned}
\end{equation*}
where the second to last inequality is a consequence of Hölder's inequality, and the last inequality follows from the Poincaré inequality.
For $(x,y)\in \supp(\pi)$ where $\pi$ is an optimal plan, we have
\begin{equation*}
    \varphi(x) = \|x-y\|^p - \varphi^c(y)
\end{equation*}
and then by differentiation with respect to $x$ we obtain
\begin{equation*}
    \nabla \varphi(x) = p (x-y) \|x-y\|^{p-2}, \qquad (x,y) \in \supp(\pi),
\end{equation*}
so that
\begin{equation*}
\begin{aligned}
\|\nabla \varphi \|_{L^{p'}(\mu)}^{p'}
= \dst\int \|\nabla \varphi \|^{p'}\ud\mu(x) = \dst\int \|\nabla \varphi \|^{p'}\ud\pi(x,y)
&= p^{p'} \dst\int \|x-y \|^{(p-1)p'}\ud\pi(x,y)\\
&= p^{p'} \dst\int \|x-y \|^{p}\ud\pi(x,y)\\
&=p^{p'} W_p^p(\mu,\nu)
\end{aligned}
\end{equation*}
and thus
\begin{equation*}
W_p^{p}(\mu,\nu)\leq \dfrac{1}{C^P_{q',p'}} \|\nabla\varphi\|_{L^{p'}(\mu)}\left[\chi^q(\nu|\mu )\right]^{\frac{1}{q}}\leq \dfrac{p}{C^P_{q',p'}}   W_p^{\frac{p}{p'}}(\mu,\nu) \left[\chi^q(\nu|\mu )\right]^{\frac{1}{q}},
\end{equation*}
which gives the first implication after dividing both sides by $W_p^{p/p'}(\mu,\nu)$.
\end{proof}

\begin{remark}\label{rmk:pis1}
    In the limit case when $p=1$ and $p'=\infty$, and $1\le q < \infty$, it is easy to prove that if all $1$-Lipschitz functions are bounded by $1/C$ then for all $\nu\in\mathcal{P}_1(\R^d)$,
    \[
    CW_1(\mu,\nu)\le \left\|\dfrac{\nu}{\mu}-1\right\|_{L^q(\mu)}.
    \]
    This follows directly from the duality of $W_1$, for $\nu\ll\mu$,
\begin{align*}
W_1(\mu,\nu)&=\sup\left\{\int f \ud\mu-\int f\ud\nu; f \text{ $1$-Lipschitz}\right\}\\
&=\sup\left\{\int f\left(\frac{\nu}{\mu}-1\right)\ud\mu; f \text{ $1$-Lipschitz}\right\}
\end{align*}
and the duality of $L^q(\mu)$
\begin{align*}
\left\|\dfrac{\nu}{\mu}-1\right\|_{L^q(\mu)}&=\sup\left\{\int f\left(\frac{\nu}{\mu}-1\right); \|f\|_{L^{q'}(\mu)}\le 1\right\}\\
\end{align*} 
Note that the result formally obtained as a limit of $\|\nabla f\|_{L^{p'}(\mu)}$ when $p'\to\infty$ is here a supremum an not an essential supremum, yielding a stronger result than the actual limit of $\|\nabla f\|_{L^{p'}(\mu)}$ suggests.
\end{remark}

\begin{remark}\label{rmk:counterT}
In the corresponding $(p,q)$-transport-divergence inequality, the constant $C^P_{q',p'}$ from the $(q',p')$-Poincaré inequality appears divided by $p$.
   This factor is optimal, as we exhibit a measure $\mu$ for which the inequality is attained up to a multiplicative constant.
   Take $\mu = e^{-x} \mathbf{1}_{\R^+}(x) \ud x$ and
   $\nu(x)= c e^{-x} \mathbf{1}_{[0,1]}(x) \ud x$ with $c =\frac{e}{e-1}$.
   We have
   \begin{equation*}
  \left[\chi^q(\nu|\mu )\right]^{\frac{1}{q}}= \Big( \frac{c-1}{c}\Big)^{\frac 1q}(1+(c-1)^{q-1})^{\frac 1q}\underset{q\to \infty}{\rightarrow}
1,
   \end{equation*}
   and
      \begin{align*}
  W_p^p(\mu,\nu)= \int_0^1|F_{\mu}^{-1}(u) -F_{\nu}^{-1}(u) |^p \ud u =\int_0^1\left|\log(1-u) -\log\!\left(1-\frac{u}{c}\right) \right|^p \ud u.
   \end{align*}

   Since $\log(1-\frac{u}{c})$ is bounded for $u\in[0,1]$, there exists $K >0$ such that

   \begin{align*}
       \int_0^1\left|\log(1-u) -\log\!\left(1-\frac{u}{c}\right) \right|^p
       \ud u &\geq  K^p \int_0^1|\log(1-u) |^p \ud u \\
       &= K^p \int_0^{+\infty} v^p e^{-v} \ud v\\
       &= K^p \Gamma(p+1) \underset{p\to \infty}{\sim} K^p p^p \sqrt{p},
   \end{align*}
   so that
         \begin{equation*}
  W_p(\mu,\nu) \underset{p\to \infty}{\sim} K p.
   \end{equation*}
   Hence
   \[
   C^T_{p,q}W_p(\mu,\nu) \underset{p\to \infty}{\sim} C^T_{p,q} K p\lesssim 1
   \]
   and thus, according to our Theorem~\ref{thm:Poincare_transport},
   \[
   O(1)\underset{p\to \infty}{=} C^P_{q',p'} K \le C^T_{p,q} K p\lesssim 1,
   \]
   showing that the factor $p$ is tight for the exponential measure $\mu$.
\end{remark}

\subsection{Transport-divergence implies Poincaré}

In this section, we show the other inequality:
\[
C^T_{p,q}\le C^P_{q',p'}.
\]

\begin{theorem}
\label{prop:Wbound}
If
 for all $\nu\in\mathcal{P}_p(\R^d)$,
    \[
  C W_p(\mu,\nu)\le \left\|\dfrac{\nu}{\mu}-1\right\|_{L^q(\mu)}
    \]
    then
    for all $f\in L^{q'}(\mu)\cap C^1(\R^d)$ such that $\int |f|^{q'-2}f\ud\mu=0$,
\[
C\|f\|_{L^{q'}(\mu)}\le\|\nabla f\|_{L^{p'}(\mu)}.
\]
\end{theorem}

\begin{remark}
    When $p=1,p'=\infty$, we can prove that the transport divergence inequality 
    \[
  C W_p(\mu,\nu)\le \left\|\dfrac{\nu}{\mu}-1\right\|_{L^q(\mu)}
    \]
    implies that all $1$-Lipschitz functions have a $L^{q'}(\mu)$-norm bounded by $1/C$, formally replacing the limit of $\|\nabla f\|_{L^{p'}(\mu)}$ as $p'\to\infty$ as a supremum and not an essential supremum.
    This is again a direct consequence of dual formulations
    \[
W_1(\mu,\nu)=\sup\left\{\int f\left(\frac{\nu}{\mu}-1\right)\ud\mu; f \text{ $1$-Lipschitz}\right\}
\]
and
    \[
    \left\|\dfrac{\nu}{\mu}-1\right\|_{L^q(\mu)}=\sup\left\{\int f\left(\frac{\nu}{\mu}-1\right); \|f\|_{L^{q'}(\mu)}\le 1\right\}.
    \]
    Indeed, suppose there exists $\bar f\in L^{q'}(\mu)$ and $1$-Lipschitz such that $\|\bar f\|_{L^{q'}(\mu)}=1/C+\eps$ for some $\eps>0$.
    We may assume, at the expense of taking a smaller positive $\eps$ that $\bar f$ is lower bounded by some $-M< 0$.
    Taking then $\nu\ll\mu$ such that $\frac{\nu}{\mu}-1=\left|\frac{\bar f}{M}\right|^{q'-2}\frac{\bar f}{M}\ge -1$ --- hence $\frac{\bar f}{M}=\left|\frac{\nu}{\mu}-1\right|^{q-2}\left(\frac{\nu}{\mu}-1\right)$, we have
    \begin{align*}
        W_1(\mu,\nu)&=\sup\left\{\int f\left(\frac{\nu}{\mu}-1\right)\ud\mu; f \text{ $1$-Lipschitz}\right\}\\
        &\ge \int \bar f\left(\frac{\nu}{\mu}-1\right)\ud\mu=M\int\left|\frac{\nu}{\mu}-1\right|^{q}\ud\mu\\
        &=M\left\|\frac{\bar f}{M}\right\|_{L^{q'}(\mu)}\left\|\frac{\nu}{\mu}-1\right\|_{L^{q}(\mu)}= \left(\frac{1}{C}+\eps\right)\left\|\frac{\nu}{\mu}-1\right\|_{L^{q}(\mu)},
    \end{align*}
    contradicting the assumption.
\end{remark}

\begin{remark}\label{rmk:counterP}
    We proved $C^T_{p,q}\le C^P_{q',p'}\le pC^T_{p,q}$ for $p,p',q,q'>1$, and that $C^P_{q',p'}\sim pC^T_{p,q}$ when $p\to\infty$ for the exponential measure (see Remark~\ref{rmk:counterT}).
    It also true that for $\mu$ the uniform distribution on $[0,1]$, $C^P_{q',p'}\sim C^T_{p,q}$ when $p\to\infty$, and thus the dependence on $p$ in the bounds is tight.
    Indeed, the uniform distribution on $[0,1]$ satisfies a $(q',p')$-Poincaré inequality for all $q'=p'$ with $p'\in[1,2]$ with constant $\sim 1$.
    Now, taking $\nu\ll\mu$, such that the density $\nu$ is smooth, then,
   
    \begin{align*}
        W_p(\mu,\nu)&=\left\|F_\nu-\Id\right\|_{L^p(\nu)},\\
        &=\left[\int_0^1\left|\int_0^t \nu(x)-1\ud x\right|^p\ud t\right]^{1/p}\\
    &\le \left\|\frac{\nu}{\mu}-1\right\|_{L^p(\mu)},
    \end{align*}
    where the inequality follows from the Jensen inequality.
    Thus, for the uniform measure on $[0,1]$, for $p\ge 2$, $C^T_{p,p}\sim C^P_{p',p'}$.
\end{remark}

The proof relies on dual representations that we first establish.
The notation and spirit of the duality is inspired from the duality \cite[Theorem~3.5]{gozlan2010transport}.

\begin{lemma}
    For $p,p'>1$, set
    \[
    \Theta(u)=-\int\inf_y\left\{\|y-x\|^p-u(y)\right\}\ud\mu(x) \quad \text{and} \quad \Upsilon(u) = \int u\ud\mu + C\left(\int |u|^{q'}\ud\mu\right)^{\frac{p'}{q'}}.
    \]
    Then, for all $\nu\in\mathcal{P}_p(\R^d)$,
    \[
    W_p^p(\mu,\nu)=\sup_{u\in C_b(\R^d)}\int u \ud \nu - \Theta(u)
    \]
    and
    \[
    \frac{\left(Cp'\right)^{1-p}}{p}\left\|\frac{\nu}{\mu}-1\right\|_{L^{q}(\mu)}^p=\sup_{u\in C_b(\R^d)}\int u \ud \nu - \Upsilon(u).
    \]
\end{lemma}
\begin{proof}
    The first equality is Kantorovich duality.
    For $p'>1$, $\Upsilon$ is strictly convex on $L^{q'}(\mu)$, so there is at most one critical point, and it is a minimum.
    Differentiating $u\mapsto \int u \ud \nu - \Upsilon(u)$ in $L^1(\mu)$ implies that a maximum is attained for $u$ such that
    \[
    \frac{\nu}{\mu}-1=Cp'\left(\int |u|^{q'}\ud\mu\right)^{\frac{p'}{q'}-1}u^{q'-1},
    \]
    Hence, since $q(q'-1)=q'$,
    \[
    \left\|\frac{\nu}{\mu}-1\right\|_{L^q}=Cp'\left(\int |u|^{q'}\ud\mu\right)^{\frac{p'}{q'p}},
    \]
    We thus get
    \begin{align*}
    \sup_{u\in C_b(\R^d)}\int u \ud \nu - \Upsilon(u)&=\int Cp'\left(\int |u|^{q'}\ud\mu\right)^{\frac{p'}{q'}-1}|u|^{q'}\ud\mu-C\left(\int |u|^{q'}\ud\mu\right)^{\frac{p'}{q'}}\\
    &=Cp'\left(\int |u|^{q'}\ud\mu\right)^{\frac{p'}{q'}}-C\left(\int |u|^{q'}\ud\mu\right)^{\frac{p'}{q'}}\\
    &=\frac{Cp'}{p}\left(\int |u|^{q'}\ud\mu\right)^{\frac{p'}{q'}}\\
    &=\frac{1}{p}\left(Cp'\right)^{1-p}\left\|\frac{\nu}{\mu}-1\right\|_{L^q}^p.
    \end{align*}
\end{proof}

\begin{lemma}\label{lem:theta_upsilon}
    For the same $p,p',q,q'$, and for all $u\in L^{q'}(\mu)\cap C^1$,
    \[
    \Theta(u)= \sup_{\nu\in\mathcal{P}_p(\R^d)} \int u \ud \nu - W_p^p(\mu,\nu).
    \]
    In addition, if $u\ge -1$, then
    \[
    \Upsilon(u)\le \sup_{\nu\in\mathcal{P}_p(\R^d)} \int u\ud \nu - \frac{\left(Cp'\right)^{1-p}}{p}\left\|\frac{\nu}{\mu}-1\right\|_{L^{q}(\mu)}^p.
    \]
\end{lemma}
\begin{proof}
We first prove the equality for $\Theta$ with two inequalities.
On the one hand, the following holds.
    \begin{align*}
        \sup_{\nu\in\mathcal{P}_p(\R^d)} \int u \ud \nu - W_p^p(\mu,\nu)&=\sup_{\nu\in\mathcal{P}_p(\R^d)} \left[\int u \ud \nu - \sup_\phi \int \phi \ud \nu + \int \phi^c \ud\mu\right]\\
        &\le - \int u^c \ud\mu,\qquad \text{taking }\phi=u,\\
        &=-\int\inf_y\left\{\|y-x\|^p-u(y)\right\}\ud\mu(x)\\
        &=\Theta(u).
    \end{align*}
    On the other hand, for $X\sim \mu$, define $Y\in\argmin\{\|y-X\|^p-u(y)\}$, and set $Y\sim \bar\nu$ and $(X,Y)\sim \pi\in\Gamma(\mu,\bar\nu)$.
    Thus,
    \begin{align*}
        \sup_{\nu\in\mathcal{P}_p(\R^d)} \int u \ud \nu - W_p^p(\mu,\nu)&= \sup_{\nu\in\mathcal{P}_p(\R^d)} \int u \ud \nu + \int u^c \ud\mu - W_p^p(\mu,\nu)-\int u^c\ud\mu\\
        &\ge \int u \ud \bar\nu + \int u^c \ud\mu - W_p^p(\mu,\bar\nu)-\int u^c\ud\mu\\
        &= \int \|x-y\|^p\ud\pi(x,y)- W_p^p(\mu,\bar\nu) - \int u^c\ud\mu\\
        &\ge - \int u^c\ud\mu\\
        &=\Theta(u),
    \end{align*}
    where the third line follows from $u(Y) + u^c(X) = \|X-Y\|^p$ by definition of $Y$, and the last inequality from $\int \|x-y\|^p\ud\pi(x,y)\ge W_p^p(\mu,\bar\nu)$.

Now, turning to $\Upsilon$, the following holds.
Suppose $\inf u>-\infty$ and set $\bar\nu\ll \mu \in\mathcal{P}_p(\R^d)$ such that
\[
    u\bar\nu=u\mu + Cp'\left(\int |u|^{q'}\ud\mu\right)^{\frac{p'}{q'}-1}|u|^{q'}\mu.
    \]
    Thus,
\begin{align*}
    \sup_{\nu\in\mathcal{P}_p(\R^d)} \int u\ud \nu - \frac{\left(Cp'\right)^{1-p}}{p}\left\|\frac{\nu}{\mu}-1\right\|_{L^{q}(\mu)}^p&\ge \int u\ud \bar\nu - \frac{\left(Cp'\right)^{1-p}}{p}\left\|\frac{\bar\nu}{\mu}-1\right\|_{L^{q}(\mu)}^p\\
    &= \int u\ud\mu + C \left(\int |u|^{q'}\ud\mu\right)^{\frac{p'}{q'}}\\
    &\geq \Upsilon(u).
\end{align*}
\end{proof}

From this, we deduce the following.
\begin{proposition}\label{prop:WtoUpsilon}
 If for all $\nu\in\mathcal{P}_p(\R^d)$,
    \[
    p^{1/p}(Cp')^{1/p'}W_p(\mu,\nu)\le \left[\chi^q(\nu|\mu )\right]^{\frac{1}{q}}.
    \]
    then, for all $u\in L^{q'}(\mu)\cap C^1(\R^d)$ with $u\ge -1$,
    \[
    \Upsilon(u)\le\Theta(u).
    \]
\end{proposition}

With the following lemma, we are able to prove the Poincaré inequality assuming $\Upsilon(u)\le\Theta(u)$.

\begin{lemma}\label{lem:upstoPoincare}
If for all $u\in L^{q'}(\mu)\cap C^1(\R^d)$, with $u\ge -1$,
    \[
\Upsilon(u)\le\Theta(u)
    \]
    then for all $u\in L^{q'}(\mu)\cap C^1(\R^d)$ bounded, such that $\int u\ud\mu=0$,
    \[
    p^{1/p}(Cp')^{1/p'}\|u\|_{L^{q'}(\mu)}\le \|\nabla u\|_{L^{p'}(\mu)}.
    \]
\end{lemma}

\begin{proof}
By approximation arguments, we can assume that $u$ is bounded.
The infimum of $y\mapsto \|y-x\|^p-u(y)$ is attained for $y^\star(x)$ such that $\nabla u(y^\star(x))=p(y^\star(x)-x)\|y^\star(x)-x\|^{p-2}$ (existence follows from coercivity of $y\mapsto\|y-x\|^p - u(y)$ since $u$ is bounded, uniqueness is not guaranteed but not necessary), and hence $\|\nabla u(y)\|^{p'}=p^{p'} \|y^\star(x)-x\|^p$.
Then, we set $T:x\mapsto y^\star(x)$ measurable.
    \begin{align*}
        \Upsilon(u)-\Theta(u)&= \int u\ud\mu + C\left(\int |u|^{q'}\ud\mu\right)^{\frac{p'}{q'}}+\int\inf_y\left\{\|y-x\|^p-u(y)\right\}\ud\mu(x)\\
        &=\int u-u\circ T\ud\mu +  C\left(\int |u|^{q'}\ud\mu\right)^{\frac{p'}{q'}}+ \frac{1}{p^{p'}}\|\nabla u\circ T\|_{L^{p'}(\mu)}^{p'}.
    \end{align*}
    Now, let 
    \[
    u(x)-u(T(x))=\langle \nabla u(\tilde T(x)),x-T(x)\rangle=-p^{-p'/p}\|\nabla u(T(x))\|^{p'-2}\langle \nabla u(\tilde T(x)),\nabla u(T(x))\rangle
    \]
    for some $\tilde T(x)\in [x,T(x)]$.
    Taking $u=\varepsilon \tilde u$ and letting $\varepsilon \to 0$, we obtain
    \begin{equation}\label{eq:ups-theta-limit}
    \frac{1}{\varepsilon^{p'}} \left[\Upsilon(\varepsilon \tilde u)-\Theta(\varepsilon \tilde u)\right] \to  C\left(\int |\tilde u|^{q'}\ud\mu\right)^{\frac{p'}{q'}}- \frac{p-1}{p^{p'}}
    \|\nabla \tilde u\|_{L^{p'}(\mu)}^{p'},
    \end{equation}
    and since $\Upsilon(u)\le \Theta(u)$ implies the left-hand side of \eqref{eq:ups-theta-limit} is non-positive, we conclude that
    \[
    C\left(\int |\tilde u|^{q'}\ud\mu\right)^{\frac{p'}{q'}} \le \frac{p-1}{p^{p'}}\|\nabla \tilde u\|_{L^{p'}(\mu)}^{p'},
    \]
    which gives the $(q',p')$-Poincaré inequality after rearranging.
\end{proof}

\begin{proof}[Proof of Theorem~\ref{prop:Wbound}]
It follows directly from Proposition~\ref{prop:WtoUpsilon} and Lemma~\ref{lem:upstoPoincare}.
\end{proof}

\section{Convex Poincaré and weak transport-divergence inequality}
\label{sec:convex}

In this section, we focus on convex Poincaré inequalities.
The natural framework associated with convex Poincaré inequalities is that of convex order and weak transport. In particular, thanks to the formulation introduced by the weak transport cost, $\overline{T}_p(\nu \mid \mu)$ can be characterized as the projection of $\mu$ onto the set of probability measures dominated by $\nu$ in the convex order.

The weak barycentric optimal transport between
$\mu$ and $\nu$
 transports each point $x$
 to a distribution whose barycenter is optimally matched to $x$.  Mass is thus rearranged by minimizing a cost depending on the averaged displacement, rather than the pointwise displacement as in classical optimal transport.

Recall the Definition~\ref{def:wbo} of weak barycentric optimal transport $\overline{\mathcal{T}}_p(\nu|\mu)$.
Note that $\overline{\mathcal{T}}_p$ is not symmetric.
As for the classical $p$-optimal transport, there is a dual formulation \cite[Theorem 2.11]{gozlan2017kantorovich}:

\begin{equation}
\overline{\mathcal{T}}_p(\nu|\mu)= \sup_{f \text{ concave}} \dst\int f^c(x)\ud\mu(x) + \dst\int f(y) \ud\nu(y),
\end{equation}
where $f^c$ is the usual convex conjugate of $f$ defined as
\begin{equation*}
f^c(x)= \inf_{y\in \R^d} \{ \|y-x\|^p -f(y)\}.
\end{equation*}

Another formulation of the weak optimal transport $\overline{T}_p(\nu \mid \mu)$ can involves the projection of $\mu$ onto the set of probability measures dominated by $\nu$ in the convex order.
Let us recall the definition of convex order.
\begin{definition}[Convex order]\label{def:convex_order}
Let $\nu,\mu\in\mathcal{P}_1(\R^d)$.
We say that $\nu$ is dominated by $\mu$ for the convex order, and we write $\nu \cx \mu$ if
\begin{equation*}
\int \varphi \, \ud\nu \leq \int \varphi \, \ud\mu
\end{equation*}
for every convex function $\varphi$.
\end{definition}
Remark that convex order is a partial order on the set of probability measures in $\mathcal{P}_1(\R^d)$.

For any $\nu\in\mathcal{P}_1(\R^d)$, we set $B_{\nu}$ as the set of probability measures that are dominated by $\nu$ in the convex order
\begin{equation*}
B_{\nu}=\{\eta \in \mathcal{P}(\R^d), \; \eta \cx \nu\}.
\end{equation*}
In \cite[Theorem 2.1]{GJ20} for $p=2$ and \cite[Lemma 6.1]{backhoff2019existence}, they prove the following useful characterization of the weak optimal transport $\overline{\mathcal{T}}_p(\nu|\mu)$: it is the distance (to the power $p$) of $\mu$ to its projection onto $B_{\nu}$, i.e.
\begin{equation}
\overline{\mathcal{T}}_p(\nu|\mu)= \inf_{\eta \in B_{\nu}} W_p^p(\eta, \mu),
\end{equation}

This perspective suggests natural candidates for the right-hand side for a transport-divergence inequality, specifically the quantities $\chi^q$ and $\tilde\chi^q$, where

\begin{equation}\label{eq:tilde_chi}
\tilde\chi^q(\nu|\mu):= \inf_{\tilde\nu \cx\nu}
\chi^q(\tilde\nu |\mu).
\end{equation}

We first state the equivalence between the two following inequalities that will prove useful in the next results.

\begin{theorem}
Let $C>0$.
For all $\nu\in\mathcal{P}(\R^d)$,
  \[
 C\left[\overline{\mathcal{T}}_p(\nu|\mu)\right]^{1/p}\le  \left[\tilde\chi^q(\nu|\mu)\right]^{1/q}
 \]
 if and only if for all $\nu\in\mathcal{P}_p(\R^d)$,
 \[
 C\left[\overline{\mathcal{T}}_p(\nu|\mu)\right]^{1/p}\le  \left[\chi^q(\nu|\mu) \right]^{1/q}.
   \]
\end{theorem}

This result is analogous to an unpublished result of Nathael Gozlan on the weak Talagrand inequality, where the KL divergence plays the role of the $\chi^q$-divergence.
\begin{proof}
It is clear that $\tilde\chi^q(\nu|\mu) \leq \chi^q(\nu|\mu)$, so that the implication from the first inequality to the second is trivial.

Now suppose that for all $\nu\in\mathcal{P}_p(\R^d)$,
\[
 C\left[\overline{\mathcal{T}}_p(\nu|\mu)\right]^{1/p}\le  \left[\chi^q(\nu|\mu) \right]^{1/q}.
\]
For $\nu\in\mathcal{P}_p(\R^d)$, take $\tilde\nu\cx\nu$ such that $\tilde\chi^q(\nu|\mu)=\chi^q(\tilde\nu|\mu)$.
Then, since $\tilde\nu\cx\nu$,
\[
\left[\overline{\mathcal{T}}_p(\nu|\mu)\right]^{1/p} \le  \left[\overline{\mathcal{T}}_p(\tilde\nu|\mu)\right]^{1/p},
\]
and thus,
\[
C\left[\overline{\mathcal{T}}_p(\nu|\mu)\right]^{1/p} \le C\left[\overline{\mathcal{T}}_p(\tilde\nu|\mu)\right]^{1/p} \le \left[\chi^q(\tilde\nu|\mu)\right]^{1/q}=\left[\tilde\chi^q(\nu|\mu)\right]^{1/q}. \qedhere
\]
\end{proof}

\subsection{Convex Poincaré and weak transport-divergence inequality}

Now, analogous to the Theorem~\ref{thm:Poincare_transport}, we can prove that measures satisfying a convex Poincaré inequality, also satisfy a weak transport divergence inequality.

\begin{theorem}
\label{thm:convex_Poincare_to_transport}
If for all $f\in L^{q'}(\mu)\cap C^1$ convex such that $\int |f|^{q'-2}f\ud\mu=0$,
\[
C\|f\|_{L^{q'}(\mu)}\le \|\nabla f\|_{L^{p'}(\mu)},
\]
then for all $\nu\in\mathcal{P}_p(\R^d)$,
    \[
 \frac{C}{p}\left[\overline{\mathcal{T}}_p(\nu|\mu)\right]^{1/p}\le  \left[\tilde\chi^q(\nu|\mu)\right]^{1/q}.
   \]
\end{theorem}

Although the proof of Theorem \ref{thm:convex_Poincare_to_transport} closely follows that of Theorem \ref{thm:Poincare_transport}, we provide the details here for completeness and for the convenience of readers who may be less familiar with weak optimal transport.

\begin{proof}
We consider $\nu \in \mathcal{P}(\R^d)$ such that $\nu$ is absolutely continuous with respect to $\mu$ and $\tilde\nu\cx\nu$. We use the duality formulation to write
\begin{equation*}
\overline{\mathcal{T}}_p(\nu|\mu)=  \dst\int \varphi^c(x)\ud\mu(x) + \dst\int \varphi(y) \ud\nu(y),
\end{equation*}
where $\varphi$ is concave and where for all $x\in \supp(\mu)$, $\varphi^c(x) \leq - \varphi(x)$, and then
\begin{equation*}
\begin{aligned}
\overline{\mathcal{T}}_p(\nu|\mu)\leq \dst\int \varphi \ud(\nu-\mu)  &=\dst\int \varphi\left(\frac{\nu}{\mu}-1 \right) \ud\mu\\
    & \leq \|\varphi\|_{L^{q'}(\mu)}\sup_{\substack{f \text{ concave} \\ \|f\|_{L^{q'}}\leq 1 }} \dst\int f \left(\frac{\nu}{\mu} -1\right) \ud\mu \\
    & \leq \|\varphi\|_{L^{q'}(\mu)}\sup_{\substack{f \text{ concave} \\ \|f\|_{L^{q'}}\leq 1 }} \dst\int  f\left(\frac{\tilde\nu}{\mu} -1\right) \ud\mu\\
    & \leq \|\varphi\|_{L^{q'}(\mu)} \left\|\frac{\tilde\nu}{\mu} -1\right\|_{L^q} \\
    &\leq \dfrac{1}{C} \|\nabla\varphi\|_{L^{p'}(\mu)}\left\|\frac{\tilde\nu}{\mu} -1\right\|_{L^q} ,
    \end{aligned}
\end{equation*}
where the last inequality follows from the convex Poincaré inequality.
Following the computation of the proof of Theorem~\ref{thm:Poincare_transport}, we get for all $\tilde \nu \cx \nu$,
\[
\dfrac{C}{p}\overline{\mathcal{T}}_p(\nu|\mu)\leq\left\|\frac{\tilde\nu}{\mu} -1\right\|_{L^q},
\]
from which we conclude by taking the infimum over $\tilde\nu\cx\nu$.
\end{proof}

\begin{theorem}
\label{thm:convex_transport_to_Poincare}
If for all  $\nu\in\mathcal{P}_p(\R^d)$,
    \[
C\left[\overline{\mathcal{T}}_p(\nu|\mu)\right]^{1/p}\le  \left[\tilde\chi^q(\nu|\mu)\right]^{1/q},
   \]
 then for all  $f\in L^{q'}(\mu)\cap C^1$ convex and such that $\int |f|^{q'-2}f\ud\mu=0$,
\[
C\|f\|_{L^{q'}(\mu)}\le\|\nabla f\|_{L^{p'}(\mu)}.
\]
\end{theorem}

\begin{proof}
Following the proof of Lemma~\ref{lem:theta_upsilon}, when $u$ is concave,
\[
\Theta(u)= \sup_{\nu\in\mathcal{P}_p(\R^d)} \int u \ud \nu - \overline{\mathcal{T}}_p(\nu|\mu),
\]
and similarly, for concave $u$,
\[
    \Upsilon(u)\le \sup_{\nu\in\mathcal{P}_p(\R^d)} \int u\ud \nu - \frac{\left(Cp'\right)^{1-p}}{p}\left\|\frac{\nu}{\mu}-1\right\|_{L^{q}(\mu)}^p.
    \]
The theorem follows the same lines as the proof of Theorem~\ref{prop:Wbound}.
\end{proof}

\subsection{Characterization of the minimizer for \texorpdfstring{$\tilde\chi^q$}{chi-q}}

The minimizer that appears in the definition of $\tilde\chi^q$ happens to have a simple characterization.

This characterization rests on the two following observations.

\begin{lemma}[Local convex-order replacement]\label{lem:replacement}
Let $\nu$ be a probability measure, and let $\theta,\theta_0$ be probability measures such that
\[
\theta_0\cx\theta.
\]
Assume that $\eps>0$ and that $\nu\ge \eps\theta$ as finite measures. Define
\[
\tilde\nu:=\nu-\eps\theta+\eps\theta_0.
\]
Then $\tilde\nu$ is a probability measure and
\[
\tilde\nu\cx\nu.
\]
\end{lemma}

\begin{proof}
The assumption $\nu\ge \eps\theta$ implies that $\tilde\nu$ is a positive measure. Since $\theta$ and $\theta_0$ are probability measures, its total mass is
\[
\tilde\nu(\R^d)=1-\eps+\eps=1.
\]
Now let $f$ be convex, with all relevant integrals well defined. Then
\[
\int f \,\ud\tilde\nu
=
\int f \,\ud\nu
-\eps\int f \,\ud\theta
+\eps\int f \,\ud\theta_0
\le
\int f \,\ud\nu,
\]
because $\theta_0\cx\theta$. This is exactly $\tilde\nu\cx\nu$.
\end{proof}

Using this Lemma~\ref{lem:replacement}, we can prove that replacing small portions of a measure by the distribution of an average of the average, creates a smaller measure in convex order.
This result is formalized in the following Proposition.

\begin{proposition}\label{prop:n-set}
Let $n\ge 2$ and $\nu$ be a probability measure, let $B_1,\dots,B_n\subset\R^d$ be disjoint Borel sets with positive $\nu$-measure, and let for all $i=1,\dots,n$,
\[
\nu_{B_i}:=\frac{\nu|_{B_i}}{\nu(B_i)},
\]
Fix $\lambda_1,\dots,\lambda_n > 0$ with $\sum_i\lambda_i=1$ and $\eps>0$ such that for any $i=1,\dots,n$,
\[
\nu(B_i)>\eps \lambda_i.
\]
Let $\nu_0$ be the distribution of $\sum_i\lambda_i X_i$ for $X_i\sim\nu_{B_i}$ independent random variables, and define
\[
\tilde\nu:=\nu-\eps\left(\sum_i\lambda_i\nu_{B_i}\right)+\eps\nu_0.
\]
Then $\tilde\nu$ is a probability measure and
\[
\tilde\nu\cx\nu.
\]
\end{proposition}

\begin{remark}
    As can be seen in the proof, the distribution $\nu_0$ can be replaced by the Wasserstein barycenter of the $\nu_{B_i}$ with weights $\lambda_i$.
\end{remark}

\begin{proof}
We can write
\[
\nu_0=\left(\sum_i\lambda_ix_i\right)\#\otimes_{i=1}^n\nu_{B_i}.
\]
For every convex function $f$,
\[
\int f \,\ud\nu_0
=
\int f\left(\sum_i\lambda_ix_i\right)\ud\otimes_{i=1}^n\nu_{B_i}(x_1,\dots,x_n)
\le
\sum_i\lambda_i\int f \,\ud\nu_{B_i}.
\]
Hence
\[
\nu_0\cx \sum_i\lambda_i \nu_{B_i}.
\]
Lemma~\ref{lem:replacement} now yields $\eta\cx\nu$.
\end{proof}

We can introduce the main tool behind our main result of the section, which states a sort of bipolar theorem for convex functions.

\begin{theorem}[Characterization of almost-surely convex functions]\label{thm:charac_conv}
    Let $\nu^\star\in\mathcal{P}_1(\R^d)$ be a probability measure.
    Let $\phi\in \mc{L}^1(\nu^\star)$ be a lower semi continuous function.
    Then, the following assertions are equivalent.
    \begin{enumerate}
        \item\label{it:asr} $\phi$ is $\nu^\star$-almost surely convex --- i.e. $\phi$ is $\nu^\star$-a.s. equal to a convex function,
        \item\label{it:cx} for all $\eta\cx\nu^\star$,
    \[
    \int \phi\ud\eta \le \int \phi \ud\nu^\star.
    \]
    \noindent Suppose in addition that $\nu^\star$ is absolutely continuous w.r.t. the Lebesgue measure.
    Then, the above assertions are also equivalent to
        \item\label{it:cxl} For all $\eta\cx\nu^\star$ such that $\eta$ is absolutely continuous w.r.t. the Lebesgue measure,
    \[
    \int \phi\ud\eta \le \int \phi \ud\nu^\star.
    \]
    \end{enumerate}
\end{theorem}

\begin{remark}\label{rmk:sphere}
    Note that in the item \ref{it:cxl} there is an added condition that the $\eta$ for which the inequality holds are absolutely continuous w.r.t. to the Lebesgue measure.
    One could ask if the restriction could be make with absolute continuity w.r.t. to $\nu^\star$.
    It is in fact not the case as the following example shows.
    Take $\nu^\star$ to be uniform on the unit sphere.
    Then, it is easy to verify to there is only one measure $\eta\cx\nu^\star$ such that $\eta\ll\nu^\star$: $\nu^\star$ itself.
    In that, any $\phi$ satisfy the inequality, but not all $\phi$ defined on the unit sphere are $\nu^\star$-almost surely convex.
\end{remark}

\begin{proof}
\textbf{\ref{it:asr} implies \ref{it:cx} and \ref{it:cxl}}

    By definition of the convex order, if $\phi$ is convex and $\eta\cx\nu^\star$, then $\int \phi\ud\eta\le\int\phi\ud\nu^\star$.
    So, we just need to prove the other implication.
    
\textbf{\ref{it:cx} implies \ref{it:asr}}

Suppose now that
\[
\int \phi \,\ud \eta \le \int \phi \,\ud \nu^\star
\]
for every probability measure $\eta$ such that $\eta\cx\nu^\star$.

Fix $n\ge 2$ and $x_1,\dots,x_n\in \supp(\nu)$ all distinct, and $\lambda_1,\dots,\lambda_n>0$ such that $\sum_i\lambda_i=1$.
Following the construction of Proposition~\ref{prop:n-set}, we define $B_i:=B(x_i,r)$ for $r>0$ small enough for the $B_i$ to be disjoint, and the corresponding measures $\nu_{B_i}$.
Choosing any admissible $\epsilon_r>0$ and applying Proposition~\ref{prop:n-set}, we obtain a measure
\[
\eta_r:=\nu^\star-\epsilon_r\left(\sum_i\lambda_i\nu_{B_i}\right)+\epsilon_r\nu_{0,r}
\]
such that
\[
\eta_r\cx\nu^\star.
\]
Hence the assumed convex-order inequality gives
\[
\int \phi \,\ud\nu_{0,r}
\le
\sum_i\lambda_i\int\phi\nu_{B_i}
\]
Using the Lebesgue differentiation theorem for probability measures, as $r\downarrow0$, the right-hand side converges $\nu^\star$-almost surely to
\[
\sum_i\lambda_i\phi(x_i)
\]
and the left hand side, $\nu^{\star\otimes n}$-almost surely converges to
\[
\lim_{r\to0}\int \phi \,\ud\nu_{0,r}=\lim_{r\to0} \frac{1}{\prod_{i=1}^n\nu(B(x_i,r))}\int_{\times_{i=1}^nB(x_i,r)} \phi\left(\sum_i\lambda_iy_i\right)\ud\otimes_{i=1}^n\nu(y_i)= \phi\left(\sum_i\lambda_ix_i\right).
\]
Hence, $\nu^{\star\otimes n}$-almost surely,
\begin{equation}\label{eq:as_conv_ineqs}
\phi\left(\sum_i\lambda_ix_i\right)\le \sum_i\lambda_i\phi(x_i).
\end{equation}

Let us deduce from these inequalities that there exists a function $\psi$, $\nu$-almost surely equal to $\phi$, that is convex.

Firstly, since \eqref{eq:as_conv_ineqs} holds for any $n\ge 1$ and any $\lambda_1,\dots,\lambda_n>0$ summing to $1$ and $\nu^{\star\otimes n}$-almost surely all $(x_1,\dots,x_n)$, there exists a set $A$ such that \eqref{eq:as_conv_ineqs} holds for any $n\ge 1$, any rational $\lambda_1,\dots,\lambda_n>0$ summing to $1$ and all $x_1,\dots,x_n\in A$, with $\nu^\star(A)=1$.

Consider the function $\phi_A := \phi \cdot \mathbf{1}_A + \infty \cdot \mathbf{1}_{A^c}$, i.e., $\phi$ restricted to $A$ and $+\infty$ outside. Define
\[
    C := \overline{\operatorname{conv}}\!\left(\operatorname{epi}(\phi_A)\right) \subset \R^n \times \R,
\]
the closure of the convex hull of the epigraph of $\phi_A$. Since $\operatorname{epi}(\phi_A)$ is closed upward in the $t$-direction, so is $\operatorname{conv}(\operatorname{epi}(\phi_A))$ and hence $C$. Therefore $C$ is itself the epigraph of a convex lower semicontinuous function, which we define as
\[
    \psi(x) := \inf\!\left\{ t \in \R : (x, t) \in C \right\},
\]
with $\psi(x) = +\infty$ if no such $t$ exists. Convexity of $\psi$ follows immediately from convexity of $C$.

We now show that $\psi=\phi$ $\nu^\star$-almost surely.

We first claim that $\psi \leq \phi$ $\nu^\star$-almost surely.
Indeed, for any $x \in A$, the point $(x, \phi(x))$ belongs to $\operatorname{epi}(\phi_A) \subset C$, so by definition $\psi(x) \leq \phi(x)$. 
Since $\nu(A) = 1$, we have $\psi \leq \phi$ $\nu$-almost surely.

Now, to prove that $\psi \geq \phi$ $\nu^\star$-almost surely, we show that for every $x \in A$ and every $(x, t) \in C$, one has $t \geq \phi(x)$, which gives $\psi(x) \geq \phi(x)$ for all $x \in A$, hence $\nu$-almost surely.

First, we claim that $C \subset \operatorname{epi}(\phi)$ when restricted to base points in $A$.
Indeed, any point $(z,s)\in C$ is the limit of points $\bigl(\sum_{i=1}^k \lambda_i x_i,\, \sum_{i=1}^k \lambda_i s_i\bigr) \in \operatorname{conv}(\operatorname{epi}(\phi_A))$ with $x_i \in A$, $s_i \geq \phi(x_i)$, and rational $\lambda_i > 0$ summing to $1$.
Then by~\eqref{eq:as_conv_ineqs} and since $\phi$ is l.s.c.,
\[
    \phi(z)\leftarrow\phi\left(\sum_{i=1}^k \lambda_i x_i\right) \leq \sum_{i=1}^k \lambda_i \phi(x_i) \leq \sum_{i=1}^k \lambda_i s_i\to s,
\]
so $(z, s) \in \operatorname{epi}(\phi)$, confirming the claim.

Therefore $\psi(x) = \inf\{t : (x,t) \in C\} \geq \phi(x)$ for all $x \in A$, i.e., $\psi \geq \phi$ $\nu^\star$-almost surely.

Hence $\phi=\psi$ $\nu^\star$-almost surely, and $\psi$ is a convex function, and the proof is complete.

\textbf{\ref{it:cxl} implies \ref{it:asr}}

The same proof holds since the measures $\eta_r$ that are used in the proof, satisfy the inequality $\int\phi\ud\eta_r\le\int\phi\ud\nu^\star$ under the extra assumptions that $\nu^\star$ is absolutely continuous w.r.t. the Lebesgue measure.
\end{proof}

We are now equipped to prove the main result of this section, which provides a necessary condition on the minimizer in the definition of the $\tilde\chi^q$ divergence.

\begin{theorem}\label{thm:chi_tilde}
    Let $q\ge 1$ and $\mu$ absolutely continuous with respect to the Lebesgue measure and with a convex support.
    Let $\nu^\star$ be a minimizer in the definition of $\tilde\chi^q(\nu|\mu)$, i.e. $\nu^\star\cx\nu$ and 
    \[
    \chi^q(\nu^\star|\mu)\le\chi^q(\eta|\mu), \quad \forall \eta\cx\nu.
    \]
    Then, $\left|\frac{\nu^\star}{\mu}-1\right|^{q-2}\left(\frac{\nu^\star}{\mu}-1\right)$ is $\nu^\star$-almost surely concave.
    \end{theorem}
\begin{remark}
    This notion is an extension to $(q-1)$-concavity usually defined for non-negative functions.
    Indeed, $f\ge 0$ is said to be $(q-1)$-concave if $f^{(q-1)}$ is concave.
\end{remark}
\begin{proof}
Take $\eta\cx\nu$ with $\eta\ll\mu$, and observe that by optimality of $\nu^\star$ and convexity of $x\mapsto |x|^q$,
\begin{align*}
\chi^q(\nu^\star|\mu)&\le \chi^q(\eta|\mu)\\
&=\int\left|\frac{\eta}{\mu}-1\right|^q\ud\mu\\
&\le\int\left|\frac{\nu^\star}{\mu}-1\right|^q\ud\mu+q\int\left|\frac{\nu^\star}{\mu}-1\right|^{q-2}\left(\frac{\nu^\star}{\mu}-1\right)\frac{\eta-\nu^\star}{\mu}\ud\mu,
\end{align*}
hence showing that for all $\eta\cx\nu$, and thus all $\eta\cx\nu^\star$ with $\eta\ll\mu$,
\[
\int \left|\frac{\nu^\star}{\mu}-1\right|^{q-2}\left(\frac{\nu^\star}{\mu}-1\right)\ud\nu^\star\le \int \left|\frac{\nu^\star}{\mu}-1\right|^{q-2}\left(\frac{\nu^\star}{\mu}-1\right)\ud\eta.
\]
Now, since this inequality holds for any $\eta\cx\nu^\star$ such that $\eta\ll\mu$, it holds for any $\eta\cx\nu^\star$ absolutely continuous w.r.t. the Lebesgue measure.
Indeed, since for any $\eta\cx\nu^\star$ with $\eta\ll\mu$, the support of $\eta$ lies in the convex hull of $\nu^\star$ which itself lies in the convex hull of the support of $\mu$ since $\nu^\star\ll\mu$.

Using the Theorem~\ref{thm:charac_conv}, we conclude that $\left|\frac{\nu^\star}{\mu}-1\right|^{q-2}\left(\frac{\nu^\star}{\mu}-1\right)$ is $\nu^\star$-almost surely concave.
\end{proof}

Note that the result is not true without the assumption the $\mu$ is absolutely continuous w.r.t. the Lebesgue measure with a convex support.
Indeed, as in the Remark~\ref{rmk:sphere}, taking any $\nu$ and $\mu$ with the unit sphere as support, then any $\eta\cx\nu$ such that $\chi^q(\eta,\mu)$ is finite is $\nu$ itself.
Hence, any $\left|\dfrac{\nu}{\mu}-1\right|^{q-2}\left(\dfrac{\nu}{\mu}-1\right)$ would be $\nu$-a.s. concave.

As a direct corollary, concavity of $\Big|\dfrac{\nu}{\mu} -1\Big|^{q-2} \left(\dfrac{\nu}{\mu} -1\right)$ characterizes that $\nu^\star=\nu$.
\begin{corollary}
Let $\mu$ absolutely continuous w.r.t. to the Lebesgue measure and with a convex support.
The two following statements are equivalent.
\begin{enumerate}
\item $\Big|\dfrac{\nu}{\mu} -1\Big|^{q-2} \left(\dfrac{\nu}{\mu} -1\right)$ is $\nu$-almost surely concave,
\item $
\tilde\chi^q(\nu|\mu) = \chi^q(\nu|\mu)$.
\end{enumerate}
\end{corollary}

\begin{proof}
The implication $2 \rightarrow 1$ is straightforward using Theorem~\ref{thm:chi_tilde}. 

We now turn to the proof of $1 \rightarrow 2$.

First, observe that by duality of $L^{q'}$,
\begin{align*}
\left[\tilde\chi^q(\nu|\mu)\right]^{1/q}
&= \left[\chi^q(\nu^\star|\mu)\right]^{1/q} \quad \text{ for } \nu^\star \cx\nu\\
&=\sup_{\substack{ \|f\|_{L^{q'}}\leq 1 }}  \int f \left(\frac{ \nu^\star}{\mu} - 1\right) \ud\mu \\
&\geq \sup\left\{ \int f \left(\frac{ \nu^\star}{\mu} - 1\right) \ud\mu; f\, \nu\text{-almost surely concave s.t. }\|f\|_{L^{q'}}\leq 1  \right\}\\
&\geq \left \|\dfrac{\nu}{\mu}-1\right\|_{L^q(\mu)}\qquad\text{ taking } f=\frac{\Big|\dfrac{\nu}{\mu} -1\Big|^{q-2} \left(\dfrac{\nu}{\mu} -1\right)}{\left\|\frac{\nu}{\mu}-1\right\|_{L^q}^{q/q'}}\\
&= \left[\chi^q(\nu|\mu)\right]^{1/q}.
\end{align*}
\end{proof}

\section{The role of \texorpdfstring{$q'$}{q'}and \texorpdfstring{$p'$}{p'} in a \texorpdfstring{($q',p'$)}{(q',p')}-Poincaré inequality on the real line}
\label{sec:role}

We now provide some intuition on the geometry of the one-dimensional space $(\R,\mu)$ when the measure $\mu$ satisfies a $(q',p')$-Poincaré inequality.
Although what we present seems part of the folklore in the community, we hope to give here clear examples to formalize that understanding.

A first observation is that the support of $\mu$ must necessarily be connected. 
Indeed, otherwise a function $f$ that is constant on the different connected components, but not constant overall, would not verify the inequality.
It is also easy to check,that if a measure $\mu$ satisfies a $(q',p')$-Poincaré inequality, applying the inequality to $|f|^s$, with $s\ge1$, for any centered and smooth $1$-Lipschitz function $f$, one has $C^P_{q',p'}\left[\int|f|^{q's}\ud\mu\right]^{1/q'}\le s\left[\int|f|^{p'(s-1)}\ud\mu\right]^{1/p'}$, showing that $\left[\int|f|^{q's}\ud\mu\right]^{1/q'}$ is finite as long as $p'(s-1)<q's$, i.e. that $f\in L^r(\mu)$ for any $r>0$ such that $\frac{1}{r}>\frac{1}{p'}-\frac{1}{q'}$, hence controlling the tail of the distribution.

To understand the geometric constraints imposed by the inequality, we investigate model situations describing both the behavior of $\mu$ at infinity and near the zeros of its density, on the real line.
Since for measures $\nu$ and $\mu$ such that the density ratio $\nu/\mu$ is bounded from above and below, then the $(q',p')$-Poincaré constants associated with $\mu$ and $\nu$ are comparable (see \cite{ledoux2018remarks}), we believe these models capture essentially when $(q',p')$-Poincaré inequalities hold.
We therefore focus on two canonical examples. 
The first one corresponds to measures whose density decays either polynomially or exponentially at infinity.
The second models measures whose density is essentially bounded away from zero everywhere on its support, except near the point 1, around which the density vanishes polynomially.

It is known since \cite{gromov1983topological} that the classical $(2,2)$-Poincaré inequality implies exponential concentration: if $f$ is a centered $1$-Lipschitz function, then
$
\mu\bigl(|f-\esp_\mu f|\ge t\bigr)
$
decays exponentially fast in $t$.
What we illustrate here is that, in dimension one, a $(q',p')$-Poincaré inequality with 
$1\leq p' \leq q'$ also implies an exponential concentration, which is a consequence of \cite[Proposition 2.5]{milmanRoleConvexityIsoperimetry2009} that links all the $(q',p')$-Poincaré inequalities.
In the case $q'<p'$ a $(q',p')$-Poincaré inequality is instead associated with heavier-tailed concentration phenomena, and the parameter $r$, defined by
\[
\frac1r=\frac1{q'}-\frac1{p'},
\]
governs the corresponding rate of decay.

We also illustrate that near a zero of the density, whatever the relative order of $p'$ and $q'$,  the polynomial behavior is constrained by the exponent $p'$: if $\mu$ satisfies a $(q',p')$-Poincaré inequality then its density can vanish near its zeros no faster than a polynomial of order $p'$.

Roughly speaking, the exponent $p'$ controls the admissible singular behavior of the density near holes in the support, while the parameter $r$ appearing in Proposition~\ref{prop:poincaré_infinity} governs the behavior of $\mu$ at infinity.

\begin{proposition}[Asymptotic behavior - case $q'<p'<\infty$ and $1\leq p'$]
\label{prop:poincaré_infinity}
Consider $\mu \in \mathcal{P}(\R)$ with the density
\begin{equation*}
g(x)=
\begin{cases}
0, \quad &x\leq 1, \\
C_{\beta} x^{-\beta}, \quad \quad  &x \geq 1,
\end{cases}
\end{equation*}
for some $\beta>1$ and $C_{\beta} =\beta-1 $.
Assume  $q'<p'<\infty$ and $1\leq p'$.
Then,  $\mu$ satisfies a $(q',p')$-Poincaré inequality if and only if $\beta \geq 1+r $, where
$$\dfrac{1}{r}= \dfrac{1}{q'} -\dfrac{1}{p'}.$$
\end{proposition}

\begin{proof}[Proof of Proposition \ref{prop:poincaré_infinity}]
 We say that a measure $\mu$ satisfies a weighted $(q',p')$-Hardy inequality if for all non-decreasing differentiable functions $h$ almost everywhere, we have
\begin{equation}
\label{eq:Hardy}
 C \|h\|_{L^{q'}(\mu)} \leq \|h'\|_{L^{p'}(\mu)}.
\end{equation}
We first show that if a measure $\mu$ satisfies a $(q',p')$-Hardy inequality, then it satisfies a $(q',p')$-Poincaré inequality. Indeed, take $f$ locally Lipschitz such that $\int f\ud\mu=0$. Define $h$ as the primitive of $|f'|$ such that $h(0)=f(0)$. Then, since $h'=|f'|$ and $|f|\leq |h|$,
\begin{equation*}
C \|f\|_{L^{q'}(\mu)} \leq C \|h\|_{L^{q'}(\mu)} \leq \|h'\|_{L^{p'}(\mu)} = \|f'\|_{L^{p'}(\mu)},
\end{equation*}
where the second inequality is the Hardy inequality for $\mu$ since $h$ is non-decreasing.
We are now ready to prove that if $\beta > r+1$, then $\mu$ satisfies the $(q',p')$-Poincaré inequality.
It is proven in
\cite{SS96} that
a measure $\mu$ satisfies a $(q',p')$-Hardy inequality on $(0,+\infty)$ if and only if
\begin{equation*}
A: =\dst\int_0^{+\infty} \left( \dst\int_0^x \mu^{1-p}(x)\ud x \right)^{\frac{r}{p}}\left( \dst\int_x^{\infty}\mu(x)\ud x\right)^{\frac{r}{p'}}  \mu(x) \ud x <\infty,
\end{equation*}
where $r$ is defined in Proposition \ref{prop:poincaré_infinity} and where $q'<p'<\infty$ and $1<p'$ and $\mu$ denote a measure and its density w.r.t. the Lebesgue measure.
A straightforward calculation shows that
\begin{equation*}
A=\dfrac{1}{(1-\beta)^{\frac{r}{p'}}}\dst\int_1^{+\infty}  x^{(1-\beta)\frac{r}{p'} }  \left( 1+ \dfrac{x^{(p-1)\beta +1}-1}{(p-1)\beta +1}\right)^{\frac{r}{p}} x^{-\beta} \ud x.
\end{equation*}
The sufficient and necessary condition for $A$ to be finite is then
$$(1-\beta) \dfrac{r}{p'} +\frac{p-1}{p} r\beta +\frac{r}{p} -\beta <-1,$$
which simplifies to $\beta > r+1$.
We have proved that if  $\beta > r+1$, then $\mu$ satisfies the $(q',p')$-Poincaré inequality.
The statement extends to  $\beta=r+1$
 by a limiting argument from the Poincaré inequality for
$\beta > r+1$.
We also obtain the case $p'=1$ by a limiting argument $p'\to 1$.
Let us now prove the reverse implication by contraposition.
Assume $\beta < r+1$, and take
\begin{equation*}
f(x)=
\begin{cases}
1, \quad &x < 1,\\
x^{\frac{\beta-1}{q'}}, \quad & x\geq 1.
\end{cases}
\end{equation*}
Then $\| f\|_{L^{q'}(\mu)}=+\infty$ and $\| f'\|_{L^{p'}(\mu)}<+\infty$, thus the $(q',p')$-Poincaré inequality fails.
\end{proof}

As the case $p'=q'=2$ already illustrates, the case $p'\leq q'$ excludes polynomial concentration and only allows faster decay rates.
Indeed, take $\mu$ as in Proposition \ref{prop:poincaré_infinity} and consider
\begin{equation*}
f(x)=
\begin{cases}
1, \quad &x < 1,\\
x^{\delta}, \quad & x\geq 1,
\end{cases}
\end{equation*}
with $\delta \in \left( \dfrac{\beta-1}{q'} , 1 +\dfrac{\beta}{p'}\right) $.
The condition $p'\leq q'$ guarantees that the interval is non empty.
Then $\| f\|_{L^{q'}(\mu)}=+\infty$ and $\| f'\|_{L^{p'}(\mu)}<+\infty$, thus the $(q',p')$-Poincaré inequality fails.

\begin{proposition}
[Asymptotic behavior - case $1\leq p'<q'<\infty$]
\label{prop:poincaré_infinity_exp}
Consider $\mu \in \mathcal{P}(\R)$ with the density
\begin{equation*}
g(x)=
\begin{cases}
0, \quad &x\leq 1, \\
\beta e^{\beta}e^{-\beta x}, \quad \quad  &x \geq 1,
\end{cases}
\end{equation*}
for some $\beta>0$.
Assume  $1\leq p'<q'<\infty$.
Then,  $\mu$ satisfies a $(q',p')$-Poincaré inequality.
\end{proposition}

\begin{proof}
It is proven in
\cite{Bradley78}  that
a measure satisfies a $(q',p')$-Hardy inequality on $(0,+\infty)$  in the context $1\leq p'<q'<\infty$ if and only if
\begin{equation*}
\sup_{y\geq 0}  \left(\dst\int_y^{+\infty} \mu(x) \ud x \right)^{\frac{1}{q'}}\left( \dst\int_0^{y}\mu^{-\frac{p}{p'}}(x)\ud x\right)^{\frac{1}{p'}}   <\infty.
\end{equation*}
A straightforward calculation shows that this condition is always satisfied for the density $\mu$ defined in Proposition \ref{prop:poincaré_infinity_exp}.
\end{proof}

\begin{proposition}[Local behavior]
\label{prop:poincaré_zero}
Consider $\mu \in \mathcal{P}(\R)$ with the density
\begin{equation*}
g(x)=
\begin{cases}
0, \quad \quad  &x \leq 0, \\
\dfrac{\alpha+1}{2}|x-1|^{\alpha}, \quad &x\in[0,2],\\
0, \quad \quad  &x \geq 2,
\end{cases}
\end{equation*}
for some $\alpha>-1$.
Assume  
$1<p'$.
Then, $\mu$ satisfies a $(q',p')$-Poincaré inequality if and only if $\alpha \leq p'-1$.
\end{proposition}

\begin{proof}[Proof of Proposition \ref{prop:poincaré_zero}]
Here again, for the case $q'<p'<\infty$, a straightforward calculation shows that if $\alpha(1-p) <-1$, i.e.\ if $\alpha > p'-1$, then $A= +\infty$, where $A$ is defined in the proof of Proposition \ref{prop:poincaré_infinity}, and otherwise
$$A= \dst\int_{-1}^{1} \left( 1+ |z|^{\alpha+1}\right)^{\frac{r}{p'}} \left( 1- |z|^{\alpha(1-p)+1}\right)^{\frac{r}{p}}  z^{\alpha} \ud z<\infty.$$
For the case $1\leq p'\leq q' <\infty$, a straightforward calculation also shows that the supremum introduced in the proof of Proposition \ref{prop:poincaré_infinity_exp} is finite is and only if $\alpha < p'-1$.
We have proved that in any case, if $\alpha < p'-1$, then $\mu$ satisfies the $(q',p')$-Poincaré inequality.
Here also, the statement extends to $\alpha = p'-1$
 by a limiting argument from the Poincaré inequality for
$\alpha < p'-1$.
Let us now prove the reverse implication by contraposition.
Assume $\alpha > p'-1$, and take
\begin{equation*}
f_{\varepsilon}(x)=
\begin{cases}
-1, \quad &x < 1-\varepsilon,\\
\dfrac{x-1}{\varepsilon}, \quad & x\in [1-\varepsilon, 1+\varepsilon],\\
1, \quad &x > 1+\varepsilon.
\end{cases}
\end{equation*}
Then $\| f_\varepsilon\|_{L^{q'}(\mu)}\underset{\varepsilon \to 0}{\rightarrow} 1$ and
$$\| f_\varepsilon'\|_{L^{p'}(\mu)}^{p'} = \varepsilon^{\alpha+1-p'}\underset{\varepsilon \to 0}{\rightarrow} 0.$$
Hence the $(q',p')$-Poincaré inequality fails when $\alpha > p'-1$.
\end{proof}

\section{Analogy with optimization problems in the Wasserstein space with \texorpdfstring{$q'=p'=2$}{q'=p'=2}}\label{sec:optim}

A point of view on functional inequalities that is becoming increasingly popular and fruitful, is to see such tools as tools of (convex) optimization on the Wasserstein space of probability measures.
We here develop somewhat systematic analogy of functional and transport inequalities with classical inequalities in optimization.

\subsection{Optimization in finite dimension}
To estimate the minimum of a functional $f$, gradient descent algorithms are the most simple algorithms and their convergence is guaranteed by
several standard optimization inequalities. Let us define the optimization inequalities in play.
Let $f:\R^d\to\R$ differentiable such that $\min f = 0$ and there exists a unique $x^\star$ such that $f(x^\star)=0$.

A function  $f$ satisfies a Polyak--Łojasiewicz \eqref{eq:PL} inequality if
\begin{equation}\label{eq:PL}\tag{PL}
     Cf(x)\le \|\nabla f(x)\|^2, \qquad  x \in\R^d.
\end{equation}

For a function $f$ (not necessarily convex), the \eqref{eq:PL} condition guarantees the convergence of the gradient descent method (see \cite{Polyak1963}).

A function $f$ satisfies a Quadratic Growth \eqref{eq:QG} inequality if
\begin{equation}\label{eq:QG}
     C\|x-x^\star\|^2\le f(x), \qquad x \in\R^d.\tag{QG}
\end{equation}

A function $f$ satisfies an Error Bound \eqref{eq:EB} inequality if
\begin{equation}\label{eq:EB}\tag{EB}
    C \|x-x^\star\|\le \| \nabla f(x)\|, \qquad x \in\R^d.
\end{equation}

A function $f$ satisfies a Weak Strong Convexity \eqref{eq:WSC} inequality if
\begin{equation}\label{eq:WSC}\tag{WSC}
     C \|x-x^\star\|^2\le f(x^\star)-f(x) - \langle  \nabla f(x),x^\star -x \rangle, \qquad x \in\R^d.
\end{equation}

In \cite{karimi2016linear}, it is proven that for a function $f$ defined on $\R^d$, we have

\begin{equation}
\label{eq:Karimi_diagram}
\eqref{eq:WSC}\implies \eqref{eq:EB}\impliedby \eqref{eq:PL}\implies \eqref{eq:QG}.
\end{equation}

\begin{remark}
The diagram presented in \cite{karimi2016linear} is slightly different since the authors use the additional assumption that $f$ has a Lipschitz gradient, whereas we do not make this assumption here as it fails to be satisfied by the two functions we present next.
\end{remark}

It is also proven in \cite{karimi2016linear} that if in addition $f$ is convex and $\nabla f$ is Lipschitz, then the reverse implications
\begin{equation*}
\eqref{eq:EB} \Longleftrightarrow \eqref{eq:PL} \Longleftrightarrow \eqref{eq:QG}
\end{equation*}
are all true.

Thanks to the Otto calculus that provides the Wasserstein space with a pseudo-Riemannian structure, optimization problems on the Wasserstein space can now be interpreted under the prism of optimization theory on $\R^d$.
For instance, the Langevin diffusion associated to a measure $\mu$ can be seen as a gradient flow of the \KL divergence $\nu\mapsto\KL(\nu|\mu)$.
The implications \eqref{eq:Karimi_diagram} also hold for functionals $F$ defined on $\mathcal{P}_2(\R^d)$, where the usual gradient is replaced by the Wasserstein gradient, and the proof is a direct adaptation of the finite-dimensional case.
In that context, the \eqref{eq:PL} inequality for a functional $F$ whose minimizer is the measure $\mu$ is
\begin{equation*}
C F(\nu) \leq \|\nabla_{W_2}F(\nu) \|_{L^2(\nu)}^2, \qquad \nu \in \mathcal{P}_2(\R^d),
\end{equation*}
and the \eqref{eq:QG} inequality is
\begin{equation*}
C W_2^2(\mu,\nu) \leq F(\nu), \qquad \nu \in \mathcal{P}_2(\R^d).
\end{equation*}
From this perspective, the tools developed for finite-dimensional optimization can be used, and convergence of Langevin diffuions can be guaranteed under some functional inequalities interpreted as optimization inequalities.

One question of interest is then to rewrite the diagram \eqref{eq:Karimi_diagram} for a specific $F$ defined on the Wasserstein space and to reinterpret the implications between the associated functional inequalities.

\subsection{The KL divergence}
When $F$ is the \KL divergence
\begin{equation*}
\KL(\nu|\mu)= \dst\int_{\R^d} \log\left( \dfrac{\nu}{\mu}\right) \ud\nu(x),
\end{equation*}
the \eqref{eq:PL} inequality is the classical log-Sobolev inequality \eqref{eq:LSI}:
\begin{equation}\label{eq:LSI}\tag{LSI}
C \; \KL(\nu|\mu) \leq \mathrm{I}(\nu|\mu),
\end{equation}
where the Fisher information $\mathrm{I}$ is defined as
\begin{equation*}\mathrm{I}(\nu|\mu)=  \dst\int \left|\nabla \log \dfrac{\nu}{\mu} \right|^2 \ud\nu(x).
\end{equation*}
The \eqref{eq:QG} inequality becomes the Talagrand transportation inequality:
\begin{equation*}
C W_2^2(\mu,\nu) \leq \KL(\nu|\mu),
\end{equation*}
and the \eqref{eq:WSC} inequality corresponds to
\begin{equation*}
C W_2^2(\mu,\nu)\leq -\KL(\nu|\mu) + \left\langle  \nabla \log \dfrac{\nu}{\mu},  T-\mathrm{Id} \right\rangle_{L^2(\nu)},
 \end{equation*}
 where $T$ is the optimal transport map from $\nu$ to $\mu$. Recalling that
  \begin{equation*}
 \|T-\mathrm{Id}\|^2_{L^2(\nu)} = W_2^2(\mu,\nu),
 \end{equation*}
 combining the \eqref{eq:WSC} inequality with the Cauchy--Schwarz inequality directly implies the HWI inequality:
\begin{equation*}
C \,\KL(\nu|\mu) \leq  W_2(\mu,\nu) \sqrt{\mathrm{I}(\nu|\mu)} - C W_2^2(\mu,\nu).
\end{equation*}

Hence, the celebrated result of \cite{OV2000} showing that the \eqref{eq:LSI} implies the Talagrand inequality can be seen as a consequence of $\eqref{eq:PL} \implies \eqref{eq:QG}$. Interestingly, the fact that for log-concave measures the Talagrand inequality implies a log-Sobolev inequality can be seen as a consequence of $\eqref{eq:QG} \implies \eqref{eq:PL}$ for convex functionals, since the \KL divergence is geodesically convex on $\mathcal{P}_2(\R^d)$ if $\mu$ is log-concave \cite[Theorem 17.15]{villani2009optimal}.

\subsection{The Chi-square divergence}
When $F$ is the $\chi^2$ divergence,
the \eqref{eq:QG} inequality is simply the $(2,2)$-transport-divergence inequality,
and the \eqref{eq:PL} inequality is the $(2,2)$-Poincaré inequality.
The latter assertion is not entirely trivial and deserves a brief explanation.
For all $\nu \in\mathcal{P}_2(\R^d)$ such that $\nu \ll \mu$, the Wasserstein gradient of $\chi^2(\cdot|\mu)$ is
\begin{equation*}
\nabla_{W_2} \chi^2(\nu|\mu)= 2 \nabla \dfrac{\nu}{\mu},
\end{equation*}
so that the \eqref{eq:PL} inequality is
\begin{equation*}
C \dst\int_{\R^d} \left( \dfrac{\nu}{\mu}-1\right)^2\ud\mu \leq  2 \dst\int_{\R^d} \left| \nabla \dfrac{\nu}{\mu} \right|^2\ud\nu, \qquad \nu\ll \mu.
\end{equation*}
By setting $f=\frac{\nu}{\mu}-1$, we obtain
\begin{equation*}
C \dst\int_{\R^d} |f|^2 \ud\mu \leq  2 \dst\int_{\R^d} \left| \nabla f \right|^2 f  \ud\mu +2  \dst\int_{\R^d} \left| \nabla f \right|^2  \ud\mu, \quad \forall f \; \text{s.t.} \dst\int f \ud\mu=0.
\end{equation*}
By replacing $f$ by $\varepsilon f$ in the previous inequality, dividing by $\varepsilon^2$ and letting $\varepsilon$ go to zero, we directly obtain the $(2,2)$-Poincaré inequality for $\mu$.

With this interpretation in mind, the fact that the $(2,2)$-Poincaré inequality implies the $(2,2)$-transport-divergence inequality follows from the implication $\eqref{eq:PL} \implies \eqref{eq:QG}$, applied to the $\chi^2$ functional.
The $\chi^2$ divergence is geodesically convex on $\mathcal{P}_2(\R^d)$ if and only if $\mu$ is log-concave \cite[Theorem 9.4.12]{ambrosio2005gradient}.
Interestingly, even without the log-concavity of $\mu$, the implication $\eqref{eq:QG} \Rightarrow \eqref{eq:PL}$ still holds in this setting: this follows from Theorem~\ref{prop:Wbound}, which shows that the $(2,2)$-transport-divergence inequality \eqref{eq:QG} implies the $(2,2)$-Poincaré inequality \eqref{eq:PL} for any measure $\mu$, without any curvature assumption. This parallels the error bound/Poincaré equivalence discussed in \cite{ledoux2018remarks}, where the $\eqref{eq:EB}$ inequality is shown to be equivalent to the Poincaré inequality for $\chi^2$.

\section{Conclusion}
We have shown that $(q',p')$-Poincar\'e inequalities are equivalent, up to a
factor $p$ that we prove to be optimal, to transport-divergence inequalities
involving the $\chi^q$ divergence, extending the classical
$(2,2)$-correspondence of
\cite{ding2015note,ledoux2018remarks,liu2020,jourdain2012equivalence} to the
whole scale of exponents $(q',p')$ introduced by Milman
\cite{milmanRoleConvexityIsoperimetry2009}. The same duality argument
extends, with essentially no extra cost, to the convex setting, where the
transport side is naturally replaced by the weak barycentric optimal
transport cost $\overline{\mathcal T}_p$ and the divergence side by the
projected divergence $\tilde\chi^q$; we further gave a geometric
characterization, via convex order, of the optimal projection defining
$\tilde\chi^q$ (Theorem~2.13). In dimension one, we
described how the exponents $q'$ and $p'$ separately control the tail
behavior of $\mu$ at infinity and the admissible degeneracy of its density
near isolated zeros, and we related the resulting hierarchy of functional
inequalities to the classical hierarchy of optimization inequalities
(Polyak--\L{}ojasiewicz, quadratic growth, error bound) on the Wasserstein
space.

\section*{Acknowledgment}
The authors would like to thank Nathael Gozlan and Paul-Marie Samson for fruitful discussions, and in particular the introduction of the quantity $\tilde\chi^q$.
This research was funded, in whole or in part, by the Agence Nationale de la Recherche (ANR), Grant
ANR-23-CE40-0017. A CC-BY public copyright license has been applied by the authors to the present
document and will be applied to all subsequent versions up to the Author Accepted Manuscript arising
from this submission, in accordance with the grant’s open access condition.

\section*{
Declaration of Use of AI
}
The overall structure of the paper was conceived without the use of AI and the paper was almost entirely written between october 2025 and March 2026.

Section 2: The proofs of Theorems 6 and 18, and of Theorems 9 and 19, were carried out without AI assistance.

Section 3: The addition of Theorem 17 was suggested by Nathael Gozlan in March 2026, based on one of his prior works.

Section 4: AI drew our attention to an equivalence between (p,q)-Poincaré inequalities and certain Hardy-type inequalities. The proofs of Propositions 28 and 29 were subsequently written without AI assistance.

Section 5: This section was written without AI and is based on unpublished notes by Thibaut Le Gouic.

Throughout the paper, AI was used to assist with English phrasing and sentence-level editing.
\printbibliography[title=References]
\end{refsection}
\end{document}